\documentclass{scrartcl}
\usepackage[T1]{fontenc}
\usepackage[latin1]{inputenc}
\usepackage{a4,amssymb,amsmath}

\newtheorem{thm}{Theorem}[section]

\newtheorem{lem}[thm]{Lemma}

\newtheorem{defn}[thm]{Definition}
\newtheorem{rem}[thm]{Remark}
\newtheorem{quest}[thm]{Question}
\newtheorem{conj}[thm]{Conjecture}

\newcommand{\norm}[1]{\left\Vert#1\right\Vert}
\newcommand{\abs}[1]{\left\vert#1\right\vert}
\newcommand{\set}[1]{\left\{#1\right\}}
\newcommand{\br}[1]{\left(#1\right)}

\newcommand{\Cp}{\mathbb C}

\newcommand{\eps}{\varepsilon}
\renewcommand{\rho}{\varrho}
\renewcommand{\phi}{\varphi}

\newcommand{\C}{\mathcal{C}}
\newcommand{\sm}{\smallskip \\}
\newcommand{\med}{\medskip \\}

\newcommand{\del}{\partial}
\newcommand{\db}{\overline\partial}

\newcommand{\tnorm}[1]{\mathrel{%
 \left\vert\kern-1pt\left\vert\kern-1pt\left\vert%
  #1%
 \right\vert\kern-1pt\right\vert\kern-1pt\right\vert}}

\newcommand{\tnorma}[1]{\mathrel{%
 \vert\kern-1pt\vert\kern-1pt\vert%
  #1%
 \vert\kern-1pt\vert\kern-1pt\vert}}

\newcommand{\tnormb}[2]{\mathrel{%
 #2\vert\kern-1pt#2\vert\kern-1pt#2\vert%
  #1%
 #2\vert\kern-1pt#2\vert\kern-1pt#2\vert}}

\def\carre{\hbox{
\vrule height 1.453ex  width 0.093ex  depth 0ex
\vrule height 1.5ex  width 1.3ex  depth -1.407ex\kern-0.1ex
\vrule height 1.453ex  width 0.093ex  depth 0ex\kern-1.35ex
\vrule height 0.093ex  width 1.3ex  depth 0ex}\,}
\def\buildo#1^#2{\mathrel{\mathop{\null#1}\limits^{#2}}}
\def\buildu#1_#2{\mathrel{\mathop{\null#1}\limits_{#2}}}

\begin{document}

\title{Fine analysis on lineally convex domains of finite type\med
{\Large\mdseries An expansion of a talk given at the Conference in honor of \\H. Skoda,
Paris, September 2005}
\footnote{MSC 2000: 32-02, 32A26, 32A35, 32F17, 32F18, 32F32, §"F45, 32T25, 32T27, 32T40, 32W05}
\footnote{Keywords: Lineally convex domains, support functions, pseudodistances, $\db$ with
nonisotropic Hölder estimates}}%
\author{Klas Diederich}%


\date{}%

\maketitle
\begin{abstract}
{\noindent\sffamily\bfseries\large Abstract}\\
A discussion of methods of nonisotropic fine quantitative complex analysis on lineally
convex domains of finite type is given. The needed support functions with best possible
estimates are considered together with the estimation of their corresponding Leray
sections with respect to nonisotropic pseudodistances. The most recent developments in this
subject are studied and open questions listed.
\end{abstract}
\section{Motivation and state of the art}

A general main topic of Complex Analysis consists in investigating the relation
between {\bfseries suitable geometric properties} of complex manifolds and fine
analytic properties of them, in particular, {\bfseries quantitative behavior of analytic
objects} on them. Of course, each time a question concerning the quantitative
behavior of an analytic object is asked, the correct notion of ''geometric
object'' has to be found which dominates its behavior.\sm
There are two main branches of this research:\sm
1) The case of {\bfseries compact manifolds}: The analytic objects here are
mostly holomorphic line bundles or, more generally, {\bfseries holomorphic
vector bundles},
and their tensor products. Certain norm conditions with respect to suitable
{\bf metrics on
the original manifold or the bundle} itself might be put on their sections. One
of the main goals of the analysis on them are {\bfseries vanishing theorems} for
analytic
cohomology, {\bf existence theorems for sections}, {\bf asymptotic behavior} for
tensor
powers. The main geometric informations needed deal with {\bf curvature conditions
(positivity resp. negativity)} or with the metrics themselves. The non-degenerate
case (concerning the metrics and their curvature) is mostly quite well understood.
All kinds of allowable {\bf degeneracies} are the main topic in recent research in this
field. The theory of multipliers and the corresponding {\bf multiplier ideal sheaves}
which originally
was introduced into Complex Analysis by J. J. Kohn in \cite{Ko4} and then carried over
to the case of complex manifolds by the work of A. Nadel (see \cite{Nad1}) is the most
important tool for this study (see \cite{Dem11}).\sm
2) The case of (relatively compact) {\bf domains in $\Cp ^n$} with more or less
smooth boundaries. Here the fundamental
goals are the existence questions for holomorphic functions with all kinds of
growth conditions and the study of their boundary behavior. The basic tool is,
of course, the study of the $\db$-Neumann problem. The necessary geometric information
sits in the complex differential geometry of the boundary, in particular, its
Levi geometry, sometimes together with topological properties of the domains and
their boundaries. \sm
Of course, there also is a mixed case (between 1) and 2)), namely the case
of (relatively compact) domains with boundaries in open
complex manifolds. Here quite new interesting phenomena occur. However, this case
is not the subject of this talk.\sm
Project 2) has been largely realized for the case of domains with non-degenerate
Levi forms. However, a lot of questions on both sides, the geometric and the analytic
aspect, remain open. In many cases the suitable geometric invariants have not yet been
found and important analytic questions are still open, since the analytic
tools are missing (see, for instance, \cite{DH1}). The situation is described in
more detail in the following subsection.

\subsection{The case of degenerate Levi forms}
\subsubsection{Geometry}
Besides the Levi geometry of the boundaries of the domains plurisubharmonic
exhaustion functions and suitable complete metrics with their curvatures are
important geometric input in this case. However, if the boundaries are sufficiently
smooth, their existence is implied by the boundary geometry. \\
Linked to the degeneracies of the Levi geometry {\bf new geometric phenomena}
appear and have to be studied:
\begin{itemize}
\item the non-diagonalizability of the Levi form.
\item The notion of type and its jumping from point to point (see \cite{An14},
\cite{An10},\cite{C4}).
\item The failure of semi-continuity of types (see \cite{An12}).
\item The jumping non-isotropical behavior of the geometry.
\end{itemize}
With respect to the degree of degeneracy there is the following scale of
severity:

\begin{enumerate}
\renewcommand{\labelenumi}{\alph{enumi})}
\item Relatively open Levi flat pieces exist in the boundary,
\item relatively open pieces of the boundary are foliated by complex manifolds of
codimension $>1$,
\item germs of complex analytic varieties of positive dimension exist in
the boundary,
\item boundary points of infinite type exist,
\item the boundary is of finite type.
\end{enumerate}
It follows from \cite{DF1} that all (relative compact) domains $D$ with smooth
$\C^\omega$-boundary are of finite type.

\subsubsection{Analysis: \mathversion{bold} $L^2$-theory.}
The analytic theory of domains with degenerate Levi form so-far is to some
extent only understood in the {\bfseries pseudoconvex case}, to which the talk restricts its
attention from now on.\sm
{\bf The $\db$-Neumann problem} is qualitatively understood in the finite type case.
For domains with $\C^\omega$-smooth boundaries its subellipticity follows
from \cite{Ko4} together with \cite{DF1}. However, the estimates for the gain
of subellipticity, so-far, are very rough and it is not even understood which
geometric invariant attached to the hypersurfaces really determines its exact
size (see \cite{DH1}). Moreover, one might ask whether the nonisotropic nature
of the geometry of these domains has been taken into account sufficiently enough
in the $\db$-Neumann problem. \\
Multipliers have been used first in \cite{Ko4}
in treating these domains. It seems to be hopeful to also use multiplier
ideal theory to get more precise information on the gain of regularity.\sm
The theory is even less developed for {\bf pseudoconvex domains with $\C^\infty$-smooth
boundary of finite type}. For them the question of subellipticity has
qualitatively
been clarified by the work of D. Catlin (see \cite{C4} and \cite{C5}). However,
the estimates on the gain of regularity are extremely rough and it seems, that
so-far good methods are missing in order to improve them for general such
domains. \sm
In addition to the finite type case which, from the point of view of
analysis just is the subelliptic case, new cases also appear on the analytic side
among weakly pseudoconvex domains with $\C^\infty$-smooth boundaries. We will
list them in the following.\sm
1. Even if subellipticity breaks down, the Neumann-operator {\bf $N$ might still
be compact}, still implying at least global regularity of the $\db$-Neumann problem.
This feature has attained a lot of attention in the last years (for a
general treatment see \cite{StrFu1}, special results are, for instance, in
\cite{Str2} and already in \cite{C6}).\sm
Local hypoellipticity at the boundary for the $\db$-problem breaks down already in
the presence of a germ of a complex analytic variety in the boundary (see \cite{DP}
and \cite{C8}).
However, the question, whether the $\db$-Neumann problem is at least {\bf globally hypoelliptic}
on all bounded weakly pseudoconvex domains with smooth boundary was much more
difficult to decide. In fact, for a long time, this was conjectured. It was
only in 1996, when M. Christ (see \cite{Chr3}) showed that for certain worm domains
as constructed in \cite{DF15} the $\db$-Neumann problem fails to be globally
hypoelliptic (see also work of D. Barrett \cite{Ba3}, \cite{Ba4} and
Chr. Kiselman \cite{K1}). It, however, should be pointed out, that it
does not become
clear from M. Christs proof (and still has not been clarified) for which
particular turning numbers the corresponding worm domain has this property. This
indicates already, how little is known about global hypoellipticity.
It becomes even more clear from the fact, that the answer to the following
question is not known:\sm
Let $D\in\in\Cp^2$ be a worm domain. Then the set
$$A:=\overline D\cap\{(z,w)\in \Cp^2:w=0\}$$
is a closed annulus and $\del D$ is strictly pseudoconvex at all other points.
Now, let $z_0 \in\del D\setminus A$ be arbitrary and let $V$ be a small neighborhood of
$z_0$. Consider any small perturbation $\tilde D$ of $D$, such that on
$\Cp^2\setminus V$ the domains $\tilde D$ and $D$ are equal. It is an open question, whether
global hypoellipticity of the $\db$-Neumann problem also necessarily
fails on $\tilde D$.
\begin{rem}
Notice, that locally and semilocally around $A$ the character of $\del D$ and
$\del\tilde D$ is the same and also the topologies agree.
\end{rem}
More generally, the following conjecture should hold true:
\begin{conj}
Let $D$ be a worm domain such that its $\db$-Neumann problem is not globally
hypoelliptic. Let $V\in\Cp^2$ be an arbitrarily small open neighbourhood of
its boundary annulus $A$ and let $\widehat{D}\subset\Cp^2$ be another
bounded smooth pseudoconvex domain such
that $\widehat{D}\cap V=D\cap V$ and such that $\del \widehat{D}\setminus V$ consists
of points of finite type only. Then the $\db$-Neumann problem on $\widehat{D}$ also
is not globally hypoelliptic.
\end{conj}

\subsubsection{Analysis: Other norms.}
For treating $\db$-problems in norms different from $L^2$, like $L^p,p\neq 2$,
Hölder norms or $\C^k$-norms, different methods are required. There are, essentially,
two approaches
\begin{enumerate}
\item A passage from $L^2$ to Hölder using techniques of microlocal
analysis (see \cite{FeKo}, \cite{FKM}).\\
\item $\db$-solving integral kernels.
\end{enumerate}
Both methods do not work on all bounded pseudoconvex domains with smooth
boundary (even not on those of finite type).\sm
For 1. the class of domains, for which this method, so-far, has been applied
successfully, is quite restricted. \\
For the construction of the integral kernels in 2. suitable Leray sections
have to be produced. They, essentially, require the existence of holomorphic
support functions depending nicely on the boundary points and satisfying good
estimates there. Unfortunately, those do not always exist, even not in $\Cp^2$ on
finite type domains, as was first shown by J.J.Kohn and L.Nirenberg in \cite{KN}
(see also \cite{Fo13}). (It has been tried by J.J.Forn{\ae}ss, \cite{Fo7} and
others to overcome this difficulty by introducing some extra techniques for the
treatment in neighborhoods of the exceptional boundary points. But also
such techniques exist so-far only in very special cases.)\sm
In this talk we want to consider a class of finite type domains which - on the
one hand - allows the construction of suitable families of support functions
and - on the other hand - is general enough to allow a large variety of degeneracies
of the Levi form, namely, lineally convex smooth domains of finite type. According
to the results of H.P.Boas, E. Straube and Yu, they do not allow failure of semicontinuity
of the type, since their general type agrees with their linear type (see
\cite{BS3}) and \cite{Y2}. However, under many other aspects, their degeneracies
can be very bad, in particular, they are not always diagonalizable.
\section{The $\db$-theory on lineally convex domains of finite type}
\subsection{Support functions}
There are essentially two different methods for the construction of smooth
families of support functions with good estimates for these domains:
\begin{enumerate}
\item A method developed by A. Cumenge in \cite{Cu3} and \cite{Cu6} based on
sharp estimates on the boundary behavior of the Bergman kernel of E. Stein
and J. McNeal (see \cite{McS}) using $L^2$-methods. (It only works in
the linearly convex case.)\\
\item A direct construction of a smooth family of holomorphic support
functions with best possible estimates (see \cite{DF19} and \cite{DF21})
\end{enumerate}
One of the essential difficulties in the construction of the support
functions lies in the following: Although they might satisfy the right
estimates in each complex tangential direction, there might be exceptional
real lines where the estimates become much worse. For the necessary estimates
of the $\db$-solving integral kernels the appearance of such exceptional real
lines is deadly. Support functions without failure of the optimal estimates
in some real directions are needed.\\
The appearance of such lines, already, becomes clear in the following simple
example:\sm
We consider in $\Cp^3=\{z=(z_1,z_2,z_3): z_j=x_j+iy_j\}$ the defining function
$$r(z):=y_1+x_2^4+x_3^6+y_3^{10}$$
and the smooth hypersurface $S=\{r=0\}$ at the point $0$. $S$ is convex, so the
real tangent space $T=\{y_1=0\}$ to $S$ at $0$ is supporting. But even in its
maximally complex subspace $T^{\Cp}$ it contains with respect to the order
of contact with $S$ two exceptional real lines. Namely, consider the complex linear
subspaces $T^{\Cp}_2:=\{z_1=z_{3}=0\}$ and $T^{\Cp}_3:=\{z_1=z_{2}=0\}$ together
with their real subspaces $T_2:\{z_1=z_3=x_2=0\}$ resp. $T_3:\{z_1=z_2=x_3=0\}$.
$T^{\Cp}_2$ has order of contact with $S$ at $0$ equal to $4$, whereas the order
of contact of $T_2$ is $\infty$, the order of contact of $T^{\Cp}_3$ is $6$,
whereas that of $T_3$ is $10$. The difficulty disappears if we replace the supporting
real hypersurface $T$ by
\begin{equation}\label{corr}
\hat{T}:= \{y_1-\eps \mathrm{Re}z_2^4 +\eps\mathrm{Re}z_3^6 =0\}
\end{equation}
for $\eps >0$ small enough. \\
Obviously, this defining function is pluriharmonic and the hypersurface defined by it has
the required order of contact with $S$ in all real directions of the tangent space.\sm
A correction analogous to (\ref{corr}) by a suitable small perturbation of the tangent
plane can be achieved at a fixed boundary point of any smooth convex hypersurface
of finite type as was shown in \cite{DMc} using an idea from \cite{DF17}.
However, such a construction does not seem to suffice for the construction of
the desired $\db$-solving integral kernels, since it does not give differentiability
of the supporting surfaces in the base point. Namely, the direction of the
exceptional real lines might jump from point to point, such that the correction terms
analogous to the ones in (\ref{corr}) would not depend differentiably on the base
point.\sm
This difficulty has been overcome by a construction by K. Diederich and J. E.
Forn{\ae}ss at first for the linearly convex case of finite type considered in
\cite{DF19} and afterwards generalized to lineally convex domains of finite
type in \cite{DF21}. The constructions are based on a new kind of analysis of
the Taylor series of convex functions which also might be interesting for different
purposes. \sm
We only state here the main result of \cite{DF21}. For this we denote by
$D=\{r(z)<0\}$ the given lineally convex domain of finite type $m$. For a point
$\zeta\in\del D$ and a vector $t\in T^{10}_{\del D}(\zeta)$ with euclidean
norm $\norm{t}=1$ we define for any $w=(w_1 ,w_2)\in\Cp ^2$
$$z_{\zeta,t}(w):=\zeta -iw_1 n_{\zeta}+w_2 t$$
where $n_{\zeta}$ is the real unit normal vector to $\del D$ at $\zeta$.\\
Furthermore, we choose a small enough open neighborhood $W_0$ of $\del D$ and put
for any point $\zeta\in W_0$
\begin{equation}
D_{\scriptstyle \zeta,t}  :={\left\lbrace {\left.w\in \Cp ^{2}:z_{\scriptstyle \zeta,t}
{\left( w\right) }  \in W_{0}: \;r_{\scriptstyle \zeta,t}
{\left( w\right) }  :=r{\left( z_{\scriptstyle \zeta,t}  {\left( w\right) }
\right) }  -r{\left( \zeta\right) }  <0\right\rbrace}  \right.}  \label{difolicon6}
  \end{equation} \sm
We put for $ j=2,\ldots ,2m$
\begin{equation}
P^{\scriptstyle j}  _{\scriptstyle \zeta,t}  {\left( w\right) }  :={\sum
_{k+l=j}  ^{}  {{{1}\over{k!}}  {{1}\over{l!}}  {{\partial ^{j}r_{\scriptstyle \zeta,t}
{\left( 0\right) }  }\over{\partial w_{2}^{k}\partial \overline{w}_{2}^{l}
}}  w_{2}^{k}\overline{w}_{2}^{l}  }}  \label{difolicon8}
\end{equation}
Notice, that the coefficients of $ P^{j}_{\zeta,t}$ are $ {\cal C}^{\infty
}$ in $ (\zeta,t)$. \smallskip \ \par
In order to be able to formulate our main result, we need the following
notation:\ \par
 \begin{defn}
\rm\label{difolicon9} \sl For any polynomial $ \sum _{j=0}^{N}\sum _{\vert
\alpha \vert +\vert \beta \vert =j}a_{\alpha \overline{\beta }}z^{\alpha }\overline{z}^{\beta
}$ on any $ \Cp ^{k}$ we put
\begin{equation}
{\left\Vert P\right\Vert }  :={\sum _{j=0}  ^{N}  {{\sum _{{\left\vert \alpha
\right\vert }  +{\left\vert \beta \right\vert }  =j}  ^{}  {{\left\vert a_{\scriptstyle
\alpha \overline{\beta }  }  \right\vert }  }}  }}  \label{difolicon49}
  \end{equation}
  \end{defn}
We then have
\begin{thm}[Di/Fo 2004]\label{SF}
$\exists$ $\hat{S}(z,\zeta)\in\C^{\infty}(\Cp ^n\times W_0)$, a holomorphic
polynomial of degree $2m$ in $z$ $\forall\;\zeta\in W_0$, such that\sm
\hspace*{7pt} i) $\hat{S}(\zeta,\zeta)=0$;\\
\hspace*{7pt} ii) For any given $\eps >0$ the function $\hat S$ can be chosen
in such a way,
that the restriction $S_{\zeta ,t}:=\hat{S}(z_{\zeta ,t}(w),\zeta )$ satisfies
\begin{equation}\label{est1}
\mathrm{Re}\, S_{\zeta, t}(w)\leq r_{\zeta ,t}(w)-\eps \sum _{j=2}^{2m}
\norm{P_{\zeta ,t}^j}\norm{w_2}^j
\end{equation}
\end{thm}
\begin{rem}
It should be stressed, that the function $\hat S$ is given by a formula which is
explicit except for the choice of two constants. Furthermore, inequality (\ref{est1})
is the best possible estimate
which can be reached on the intersection of $D$ with all the $\Cp ^2$'s as spanned
by all $n_{\zeta}$ and $t$.
\end{rem}

\subsection{The pseudodistance}
\subsubsection{The definition of the pseudodistance}
Following ideas of E. Stein and from \cite{BNW}, J. McNeal introduced in \cite{Mc2}
for linearly convex domains of finite type a pseudometric which reflects very precisely
the non-isotropic nature of the geometry of these domains. It has been generalized
to lineally convex domains of finite type in \cite{Conr1} (using ideas from
\cite{Hef1}). It is defined in the following way:\sm
Let $\rho$ be a defining function of the lineally convex domain
$D\subset\subset\Cp ^n$ of finite
type. Put for any $\eps >0$, any point $\zeta$ close enough to $\del D$ and any
vector $\gamma\in\Cp ^n$
$$\tau (\zeta ,\gamma ,\eps):=\max\set{c:\abs{\rho (\zeta +\lambda\gamma)-\rho (\zeta )}
<\eps\;\forall\lambda\in\Cp :\abs{\lambda}<c}$$
With this we next choose what we call an {\bf $\eps$-extremal basis} of $\mathbb{C}^n$
at $\zeta$ in the following way:\sm
We choose as $v_1 (\zeta ,\eps)$ the unit vector orthogonal to the level set of
$\rho$ passing through $\zeta$. Then we restrict attention to the linear subspace
$H_1$ orthogonal to $v_1 (\zeta ,\eps)$ and in it we choose a unit vector $v_2$
pointing in a direction $\gamma\in H_1$ such that $\rho (\zeta +\gamma)=\pm\eps$
and $\tau (\zeta ,\gamma ,\eps)$ is maximal among such $\gamma$. This procedure is
repeated until the orthonormal basis $\br{v_1(\zeta ,\eps),\ldots ,v_n (\zeta ,\eps)}$
is complete.
\begin{rem}
Notice, that, in general, the dependence of the $\tau (\zeta ,\gamma ,\eps)$ and
the $\eps$-extremal basis on the point $\zeta$ is not even continuous.
\end{rem}
The {\bf $\mathbf \eps$-distinguished polydiscs} and their versions scaled by
$A>0$ are defined as
$$AP_{\eps}(\zeta):=\set{z=\zeta +\sum\lambda _kv_k(\zeta ,\eps):\abs{\lambda_k}
\leq A\tau_k(\zeta ,\eps)\mbox{ for }k=1\ldots ,n}$$
Finally, the non-isotropic {\bf pseudodistance} is defined as
$$d(\zeta ,z):=\inf\set{\eps :z\in P_{\eps}(\zeta)}$$
\subsubsection{Properties of the pseudodistance}
Although the pseudodistance is not continuous, it can be shown to satisfy certain
uniform estimates making it, nevertheless, possible to do analysis with it. One has:
\begin{enumerate}
\renewcommand{\labelenumi}{\roman{enumi})}
\item $\exists\: c>0:cP_{\abs{\rho (\zeta)}}\subset D\:\forall\;\zeta$
\item If $\gamma =\sum_1^n a_j v_j(\zeta ,\eps)$ $\Rightarrow$
$$\frac{1}{\tau (\zeta ,\gamma ,\eps)}\approx \sum_{j=1}^n
\frac{\abs{a_j}}{\tau _j(\zeta ,\eps)}$$
\item $\forall\: k>0$ $\exists$ constants $c(k),C(k)$ such that
\begin{equation}\label{engulfing}
c_kP_{\eps}(\zeta)\subset P_{k\eps}(\zeta)\subset C(k)P_{\eps}(\zeta)
\end{equation}
\item $\forall\; z\in P_{\eps}(\zeta )\:\Rightarrow\tau (\zeta ,\gamma ,\eps )
\approx \tau (z ,\gamma ,\eps)$
\item \begin{eqnarray*}
d(z,\zeta )& \approx & d(\zeta ,z)\\
d(z,\zeta )&\lesssim & d(z,w)+d(w,\zeta )
\end{eqnarray*}
\end{enumerate}
Here the relation $\approx$ respectively $\lesssim$ stand for the corresponding
strict relations $=$ resp. $\leqq$ up to constants uniform in the choice of $\zeta$
chosen from a small enough neighborhood of the boundary $\del D$ and $\eps >0$.
\subsubsection{Estimates relative to the pseudodistance}
The pseudodistance introduced above reflects exactly the non-isotropic geometry of
the corresponding domain $D$. It, therefore, is not astonishing that the estimates
required for the quantitative solutions of $\db$, become simple and easy to use -
just as the estimates with respect to the euclidean metric in the strictly
pseudoconvex case. We will show this by giving a few examples.
\begin{lem}\label{ES}
Define $\zeta\in D\cap W$ and $\eps >0$ the set $P_{eps}^0(\zeta):=
CP_{\eps}(\zeta)\setminus \frac{1}{2}P_{\eps}(\zeta)$ ($C=C(1)$ being the constant
of (\ref{engulfing})) and let $\pi$ be the orthogonal projection to $\del D$.
Then one has $\forall\: z\in D\cap U$
\begin{eqnarray*}
\abs{S(z,\zeta )}&\gtrsim &\eps\:\forall\:\zeta\in\del D\cap P^0_{\eps }(\pi(z))\\
\abs{S(z,\zeta )}&\gtrsim &\abs{\rho (z)}\:\forall\:\zeta\in\del D\cap P_{\abs{\rho (z)}}(\pi(z))
\end{eqnarray*}
\end{lem}
From the explicit formula for the family of support functions $S(z,\zeta )$ one easily
gets a corresponding {\bf Leray decomposition},
for which one has
\begin{equation}
S(z,\zeta )=\sum_{j=1}^nQ_j(z,\zeta )(z_j-\zeta _j)
\end{equation}
(for details see \cite{DF21} and \cite{DFi2}). In order to get the estimates required
for the corresponding $\db$-solving Cauchy-Fantappié kernels, derivatives of the
Leray section $(Q_1(z,\zeta ),\ldots ,Q_n(z,\zeta ))$ have to be estimated. Using
the pseudometric in a consequent way, this calculus becomes very natural:\\
For a point $\zeta _0\in\del D$ and arbitrary $z,\zeta \in P_{\eps}(\zeta _0)$ we
choose a linear transformation $\Phi$ giving the $\eps$-extremal coordinates at
$\zeta _0$ and put
\begin{equation*}
w:=\Phi (z-\zeta _0)\text{ and }\eta :=\Phi (\zeta -\zeta _0)
\end{equation*}
($\Rightarrow$ $\abs{\eta_k}\leq C\tau _k(\zeta _0,\eps );\:\abs{w_1}\leq C;
\abs{w_k}\leq C\tau _k(\zeta _0,\eps )$)\sm
We denote the Leray section in the new coordinates by $Q^*_1(w,\eta ),\ldots ,
Q^*_n(w,\eta )$ and get
\begin{lem}\label{EQ}
\begin{eqnarray*}
\abs{Q^*_k(w,\eta )}&\lesssim &\frac{\eps }{\tau _k(\zeta _0,\eps)}\\
\abs{\frac{\del}{\del w_i}Q^*_k(w,\eta }&\lesssim &
\frac{\eps }{\tau _k(\zeta _0,\eps)\tau _i(\zeta _0,\eps)}\\
\abs{\frac{\del }{\del\eta _j}Q^*_k(w,\eta }&\lesssim &
\frac{\eps }{\tau _k(\zeta _0,\eps)\tau _j(\zeta _0,\eps)}\\
\abs{\frac{\del ^2}{\del w_i\del\eta _j}Q^*_k(w,\eta }&\lesssim &
\frac{\eps }{\tau _k(\zeta _0,\eps)\tau _j(\zeta _0,\eps)\tau _i(\zeta _0,\eps)}
\end{eqnarray*}
\end{lem}

\subsection{A first result on solving \mathversion{bold}$\db$ with nonisotropic estimates}
At first, using the above-mentioned smooth family of holomorphic support functions,
quantitative results on solving $\db$ where given only with respect to isotropic Hölder
norms. We mention as examples \cite{DFF, Hef2, Cu3, Cu6}. However, it is natural,
that also the norms measuring the solutions to the $\db$-equation have to respect the
nonisotropic nature of the geometry. As a typical example we define new Hölder
norms by
\begin{defn}
For any $\mu >0$ and small $\eps >0$ we define
\begin{equation}\nonumber
\tilde{\Lambda}^{\mu ,\eps}(D):=\set{h\in C^0(D):\abs{h(z_0)-h(z_1)}\leq
C_h\max\set{d(z_0,z_1)^{\mu},\abs{z_0-z_1}^{1-\eps }}}
\end{equation}
\end{defn}
\begin{rem}
For any given $h$ the smallest constant $C_h$ possible in this inequality is called
the corresponding Hölder norm $\norm{h}_{\mu ,\eps}$. The term $\abs{z_0-z_1}^{1-\eps }$
only appears in this definition only for technical reasons, namely, in order to take
care of the case, when the points $z_0$ and $z_1$ are far apart. For the
understanding of the results, it suffices to neglect this term and, hence, the role
of $\eps$ in the definition of the Hölder norm.
\end{rem}
One has with respect to these nonisotropic Hölder norms (see also \cite{CuFr} for
a slightly weaker result)
\begin{thm}[Di/Fischer 2004]\label{ni}
For any $1\leq q\leq n$ $\exists$ a continuous linear operator
$T_q:L^\infty _{0,q}\rightarrow \tilde{\Lambda}^{\frac{1}{m},\eps }_{(0-q-1)}$
such that $\db T_qf=f\:\forall f:\db f=0$.
\end{thm}
\begin{rem}
a) The result of the theorem only is a test case for a whole series of possible
similar results. For instance, correct nonisotropic Hölder norms
have to be defined for $\db$-closed $(0,q)$-forms and the corresponding $(0,q-1)$-forms solving $\db$ have to be found
solving them (see also below).\\
b) The proof of Theorem \ref{ni} is quite natural once one has the support functions
of Theorem \ref{SF}, the pseudodistance and the estimates of Lemmas \ref{ES}
and \ref{EQ}. In fact, it seems, that it might be possible to formalize the theory
by proving a technical general proposition which says, that, on a given domain, one always has
such such a theorem if for it support functions and a pseudodistance
exist such that the corresponding estimates of the Lemmas \ref{ES} and \ref{EQ} hold
(i.e. other special properties of the domain do not enter into the proof).
\end{rem}
\subsection{\mathversion{bold} $\db$ and differentiability up to the boundary}
A general result of J. J. Kohn (see \cite{Ko2}) guarantees that on any weakly
pseudoconvex domain with smooth boundary any $\db$-closed $(0,1)$-form which is
$\C^{\infty}$ up to the boundary can be solved by a function $\C^{\infty}$ up
to the boundary. However, there is no good $\C^k$-estimate for the solution.\sm
For strictly pseudoconvex smooth domains this problem has first been solved
by G. Henkin (see \cite{He5}). Recently, W. Alexandre (see \cite{Ale1, Ale2}) showed a first
result in this direction on smoothly bounded, bounded convex domains of finite
type. He uses as a main tool again the above mentioned smooth family of
holomorphic support functions.
\begin{thm}
Let $D\subset\subset \Cp ^n$ be a linearly convex domain with $\C^{\infty}$-smooth
boundary of finite type $\leq m$. Then there is for any $1\leq q\leq n$ and
any $k=0,1,2,\ldots$ a bounded linear operator
\begin{equation*}
T_q:\C^k_{(0,q)}(\overline D)\cap\ker\db\rightarrow\C^{k,\frac{1}{m}}_{(0,q-1)}(\overline D)
\end{equation*}
such that
$$\db T_qf=f$$
\end{thm}
\begin{rem}
\hspace{7pt}a) In the Theorem the spaces $\C^k_{(0,r)}(\overline D)$ are provided with the $\sup$-norm
on $\overline D$ over all partial derivatives up to the total order $k$ of all
coefficients and
the norm on the space $\C^{k,\frac{1}{m}}_{(0,q)}(\overline D)$ is the sum of the
norm on
$\C^{k}_{(0,q)}(\overline D)$ and the isotropic Hölder norms of order $\frac{1}{m}$
on $\overline D$ of all partial derivatives of order $k$ of all the coefficients
of the forms.\\
\hspace{7pt}b) Recently, K. Di. and B. Fischer have carried over this result to
lineally convex domains of finite type and, more importantly, to the analogous
nonisotropic Hölder norms on $\mathrm{Im}\,T_q$. The corresponding preprint will
appear soon and will contain some further research on the subject.
\end{rem}

\subsection{The extension problem}
In \cite{Maz2} E. Mazzilli gave an example of a pseudoellipsoid $D$ together with
a complex-analytic algebraic variety $X$, smooth in an open neighborhood of
$\overline D$ and intersecting $\del D$ transversally, such that, nevertheless,
there is a bounded holomorphic function $f$ on $X\cap D$ which does not
extend to a bounded holomorphic function on $D$. This strange
phenomenon was, afterwards, more closely studied in \cite{DiMaz1}. Later,
K. Di and E. Mazzilli showed in \cite{DiMaz4}, that, in contrast to this,
every bounded holomorphic function on $X\cap D$ extends to a bounded
holomorphic function if $X$ is affine linear and if $D$ is just a bounded linearly
convex domain with smooth boundary of finite type. \sm
The question which submanifolds allow a bounded holomorphic extension of all
bounded holomorphic functions has first been studied by E. Mazzilli. He
gave a sufficient condition using different type notions for the
intersection $\del D\cap X$. However, his condition is too strong to be also
necessary. Recently, W. Alexandre used in \cite{Ale3} the nonisotropic pseudometric
as introduced above to give a sufficient condition on $X\cap D$ for bounded
extendibility (for linearly convex domains of finite type) with respect to
the $\sup$-norm which is very close to also being necessary. Since the formulation
of the condition is rather technical we refer the interested reader to the original
paper.

\subsection{Some open problems}
Besides the work mentioned so-far, many more questions of quantitative complex
analysis have been studied on linearly (or even lineally) convex domains
of finite type. In addition to the above, we mention the following selection:
\begin{itemize}
\item Boundary behavior of $H^p$-functions, see \cite{DIBIF1}.
\item Bounded solvability of $\db$ with respect to $L^p$-norms with $p\neq 2$, see
\cite{Fi}.
\item Characterization of zero sets of the Nevanlinna class, see
\cite{BrChDu, Cu4, DiMaz2}
\item Dependence of isotropic Hölder and $L^p$-norms on Catlins multitype,
see \cite{Hef1}.
\item $(0,1)$-forms with finite Bruna-Charpentier-Dupain-norm can be
solved in $L^1(\del D)$, see \cite{BrChDu, Cu4, DiMaz2}.
\end{itemize}
As closing remarks of this survey we mention several important open questions
of the subject:
\begin{quest}
\begin{enumerate}
\item\label{IP} Generalize all results known only for linearly convex domains of finite type
to lineally convex domains of finite type.
\item Investigate for $\mu>1$ the relation between the spaces $\C^{l,\mu}(\overline D)$ and
and $\C^{l+k(\mu),\mu-k}(\overline D)$, where $k(\mu)$ is the largest integer $\leq\mu$ and
all Hölder norms are taken with respect to the pseudodistance on $D$.
\item Solve the $\db$-equation with best possible estimates for forms with coefficients
in $\widetilde{\Lambda}^{\mu}$ respectively in $\C^{k,\mu}$.
\item Study more general classes of domains to which the construction of a smooth
family of holomorphic support functions satisfying best possible estimates
of Di-Forn{\ae}ss can be carried over. \label{class}
\end{enumerate}
\end{quest}
\textit{Added in proof:} Recently, after finishing this survey, the author has been informed, that
a preprint of J. Michel appeared,
in which he generalizes the Diederich-Forn{\ae}ss construction of support functions to so-called
K-convex domains.

------------------------------------------------------------------

\makeatletter \renewcommand{\@biblabel}[1]{\hfill#1.}\makeatother
\newcommand{\bysame}{\leavevmode\hbox to3em{\hrulefill}\,}

\ \par
\renewcommand{\baselinestretch}{0.7} \Large \normalsize \ \par
  \vspace{-\baselineskip}\begin{tabbing}
 \=\mit \footnotesize Klas\ Diederich\ \ \ \ \ \ \ \ \ \ \ \ \ \ \ \ \ \ \ \ \ \ \ \ \ \ \ \ \ \ \ \ \ \ \ \ \ \ \ \ \ \ \ \ \ \ \ \  \=\mit
\footnotesize \\\mit
 \>\mit \footnotesize Mathematik  \>\mit \footnotesize \\\mit
 \>\mit \footnotesize Universit\"at Wuppertal  \>\mit \footnotesize \\\mit
 \>\mit \footnotesize Gausstr. 20  \>\mit \footnotesize \\\mit
 \>\mit \footnotesize D-42097 Wuppertal \>\mit \footnotesize \\\mit
 \>\mit \footnotesize Germany  \>\mit  \\\mit
  \end{tabbing}\vspace{-2\baselineskip}\mit
\ \par
\end{document}